\documentclass[12pt]{amsart}
\usepackage{amsmath,amsthm,amssymb,bbm,dsfont}
\usepackage{graphicx}
\usepackage{subfig}
\usepackage{cite,float}

\numberwithin{equation}{section}
\newtheorem{theorem}{Theorem}[section]
\newtheorem{lemma}[theorem]{Lemma}

\theoremstyle{remark}
\newtheorem{remark}[theorem]{Remark}
\newtheorem{example}[theorem]{Example}
\newtheorem{definition}[theorem]{Definition}

\makeatletter
\@namedef{subjclassname@2020}{%
  \textup{2020} Mathematics Subject Classification}
\makeatother
\begin{document}
\sloppy
\title{Squeezing Function on Infinitely Connected Planar Domain}

\author[A. Kumar]{Akhil Kumar}
\address{Department of Mathematics, University of Delhi,
Delhi--110 078, India}
\email{98.kumarakhil@gmail.com}

\begin{abstract}
We provide explicit expression of squeezing function for infinitely connected planar domain obtained by removing a convergent sequence of points from the unit disk converging to the boundary of unit disk. We also discuss Fridman invariant associated with this domain as well as some examples of squezzing functions corresponding to polydisk.  
\end{abstract}
\maketitle

\section{Introduction}
In 2012, Deng-Guan-Zhang in their paper \cite{2012}, introduced the notion of squeezing function by following the work of Liu-Sun-Yau \cite{Yau2004}, \cite{Yau2005} and Yeung \cite{yeung}. One of the many interesting problems is to find explicit formulae for the squeezing functions associated with different domains in one and more complex dimensions. Recently Ng-Tang-Tsai \cite{recent} and Gumenyuk-Roth \cite{roth} have written articles to this effect. In one dimension the case of infinitely connected domain is largely unexplored as far as the explicit formulae is concerned. The main result of this paper is
\begin{theorem}\rm\label{main1}
Let  $\Omega = \mathbb{D}\setminus A$ be a infinitely connected domain, where $A = \lbrace a_{n}:  n\in \mathbb{N} \rbrace$ is a sequence in $\mathbb{D}$ converging to $\partial\mathbb{D}$(boundary). Then, the squeezing function on $\Omega$, is given by
$$S_{\Omega}(z) = \inf_{w\in A} \left\lbrace \left|\dfrac{w-z}{1-\overline{z}w}\right|\right\rbrace.$$
\end{theorem}  
Here $S_{\Omega}(z)$ is expression used for squeezing function whose definition in general form is given as:
\begin{definition}
Let $\Omega\subseteq \mathbb{C}^n$ be a bounded domain. For $z\in{\Omega}$ and a holomorphic embedding $f:\Omega\to \mathbb{B}^n$  with $f(z)=0$, define 
$$S_{\Omega}(z,f):=\sup \{r:\mathbb{B}^n(0,r)\subseteq f(\Omega)\},$$
where $\mathbb{B}^n$ denotes unit ball in $\mathbb{C}^n$ and $\mathbb{B}^n(0,r)$ denotes the ball centered at origin with radius $r$ in $\mathbb{C}^n$.
The squeezing function on $\Omega$, denoted by $S_{\Omega}$,
is defined as
$$S_{\Omega}(z):=\sup_f\{S_{\Omega}(z,f)\},$$ 
where supremum is taken over all holomorphic embeddings $f:\Omega\to \mathbb{B}^n$ with $f(z)=0.$
\end{definition}

The analogous theory of squeezing function when holomorphic embeddings are taken into unit polydisk instead of unit ball, has been given by Gupta and Pant \cite{polydisc}.  The definition takes the form as:
\begin{definition} 
Let $\Omega\subseteq \mathbb{C}^n$ be a bounded domain. For $z\in{\Omega}$ and a holomorphic embedding $f:\Omega\to \mathbb{D}^n$  with $f(z)=0$, define 
$$T_{\Omega}(z,f):=\sup \{r:\mathbb{D}^n(0,r)\subseteq f(\Omega)\},$$
where $\mathbb{D}^n$ denotes unit polydisk and $\mathbb{D}^n(0,r)$ denotes the polydisk centered at origin with radius $r$ in $\mathbb{C}^n$.
The squeezing function corresponding to polydisk on $\Omega$, denoted by $T_{\Omega}$,
is defined as
$$T_{\Omega}(z):=\sup_f\{T_{\Omega}(z,f)\},$$ 
where supremum is taken over all holomorphic embeddings $f:\Omega\to \mathbb{D}^n$ with $f(z)=0.$
\end{definition}

It is clear from the definitions that for any bounded domain $\Omega$, squeezing functions are biholomorphic invariant.  It is interesting to note that the domain considered in our main theorem is not a  holomorphic homogenous regular domain. A domain is said to be holomorphic homogenous regular(HHR) domain if squeezing function on it admits positive lower bounds. For other important properties of squeezing functions one can refer the following papers: \cite{deng2019}, \cite{uniform-squeezing},\cite{kaushalprachi}, \cite{recent}, \cite{rong}, \cite{rong2}, \cite{rong3} and references therein.

Recently Rong and Yang \cite{rong} have generalized squeezing function by taking holomorphic embeddings into a bounded, balanced and convex domain in $\mathbb{C}^{n}$. Recently Gupta and Pant \cite{d-balanced}, introduced more general notion of $d$-balanced squeezing function by taking the embeddings into bounded, $d$-balanced, convex domain. Now the question comes, what is the importance of the object $T_{\Omega}$, when there are more general set ups available. It seems to us that the object $T_{\Omega}$ is more amenable to the computation in comparison to the object associated with arbitrary balanced domains. Our next result talks about Fridman invariant.

\begin{theorem}\rm\label{main2}
For  domain $\Omega = \mathbb{D}\setminus A$, where $A$ as in Theorem 1.1, the Fridman invariant on $\Omega$, is given by
$$h_{\Omega}^{c}(z) = \inf_{w\in A} \left\lbrace \left|\dfrac{w-z}{1-\overline{z}w}\right|\right\rbrace.$$
\end{theorem}

In \cite{frid1979, frid1983}, Fridman introduced a holomorphic invariant known as Fridman invariant. For domain $\Omega\subseteq \mathbb{C}^n$ and $z\in{\Omega}$, Fridman invariant denoted by
$h_{\Omega}^{d}$, is defined as
$$h_{\Omega}^{d}(z):= \inf\left\lbrace \dfrac{1}{r} : B_{\Omega}^{d}(z, r)\subseteq f(\mathbb{B}^{n}), f: \mathbb{B}^{n}\to \Omega \right\rbrace,$$
where $f: \mathbb{B}^{n}\to \Omega$ is a holomorphic embedding and $B_{\Omega}^{d}(z, r)$ is a ball centered at $z$ with radius $r$ with respect to Carathéodory/Kobayashi pseudometric.  For comparison purpose, we will follow Nikolov and Verma \cite{kaushal}, by considering Fridman invariant replacing infimum by supremum and $\dfrac{1}{r}$ by $\tanh r$(both in their uses lead to the same thing). For Carathéodory and Kobayashi pseudometric, we denote Fridman invariants by $h_{\Omega}^{c}$ and $h_{\Omega}^{k}$ respectively. Nikolov and Verma have established a relation between both these invariant functions as $S_{\Omega}(z)\leq h_{\Omega}^{c}(z)\leq h_{\Omega}^{k}(z)$, $z\in \Omega$. It is natural to ask, when does equality hold. Note that if $\Omega$ is holomorphic equivalent to $\mathbb{B}^{n}$, then equality holds trivially. In \cite{rong2}, Rong and Yang have shown that for bounded, balanced, convex and homogeneous domain $\Omega$ in $\mathbb{C}^{n}$, $S_{\Omega}(z) = h_{\Omega}^{c}(z) = h_{\Omega}^{k}(z)$. In this paper, we see that $S_{\Omega}(z) = h_{\Omega}^{c}(z)$ for $\Omega = \mathbb{D}\setminus A$.   

\section{Proof of theorems} 
Before stating we would like to fix some notations. Let $\mathbb{D}$ be the unit disk. Then for a domain $\Omega\subseteq \mathbb{C}^{n}$ and $z_{1}, z_{2}\in \Omega$, the Carathéodory pseudodistance on $\Omega$ is defined as $c_{\Omega}(z_{1}, z_{2}) = \sup\lbrace \tanh^{-1}|\alpha| : f:\Omega\to \mathbb{D}\  \ holomorphic, 
f(z_{1})=0, f(z_{2})= \alpha \rbrace$. Note that, for $z = (z_{1}, \ldots, z_{n})$ and $w = (w_{1}, \ldots, w_{n})$, we have
$$\tanh c_{\mathbb{D}^{n}}(z, w) =  \max_{1\leq i \leq n}\left| \dfrac{w_{i} - z_{i}}{1-\overline{z_{i}}w_{i}}\right|.$$
Let $d_{\rho}(z_{1}, z_{2})$ denote the Poincaré distance between $z_{1}$ and $z_{2}$ in $\mathbb{D}$. We consider a real valued function  $\sigma:[0, 1)\to \mathbb{R}_{\geq 0}$ defined as
$\sigma(x) = \dfrac{1}{2}\log\dfrac{1+x}{1-x}$. It is clear that $\sigma$ is strictly increasing function and for $z\in\mathbb{D}$, $\sigma(|z|)$ is the Poincaré distance between 0 and $z$. One can observe that the inverse of this function can be expressed as $\tanh(x)$.

\begin{proof}[\underline{Proof of Theorem \rm{\ref{main1}}}]
It is clear that the infimum exists and we need to show that it can not be 0. Let $\inf_{w\in A} \left\lbrace \left|\dfrac{w-z}{1-\overline{z}w}\right|\right\rbrace$ be 0, then there exists a sequence $\lbrace a_{k} \rbrace$ in $A$ such that $\left\lbrace \left|\dfrac{a_{k}-z}{1-\overline{z}a_{k}}\right|\right\rbrace$ converges to 0. Then for $\epsilon > 0$, there is some $N\in \mathbb{N}$ such that for all $k\geq N$, we have
$$\dfrac{|a_{k}-z|}{2} \leq\left|\dfrac{a_{k}-z}{1-\overline{z}a_{k}}\right| < \dfrac{\epsilon}{2}.$$
This implies $\lbrace a_{k} \rbrace$ converges to $z$ which is not possible. Now for $z\in \Omega $, consider the function $f :\Omega \to \mathbb{D}$ 
$$f(\zeta) = \dfrac{\zeta-z}{1-\overline{z}\zeta}.$$
This is a holomorphic embedding with $f(z) = 0$. As we can see that the disk centered at zero with radius $\inf_{w\in A} \left\lbrace \left|\dfrac{w-z}{1-\overline{z}w}\right|\right\rbrace$ contained in $f(\Omega)$, this implies
\begin{align}
S_{\Omega}(z) \geq \inf_{w\in A} \left\lbrace \left|\dfrac{w-z}{1-\overline{z}w}\right|\right\rbrace.
\end{align}
For other inequality, let $S_{\Omega}(z) = c$ for some $c\in (0 ,1]$, then by {\cite[Theorem~2.1]{2012}}, there is a holomorphic embedding $f :\Omega \to \mathbb{D}$ with $f(z) = 0$ such that $\mathbb{D}(0 , c)\subseteq f(\Omega)$. Let $\Omega_{1}= \Omega \cup \lbrace a_{1} \rbrace$, $\Omega_{2}= \Omega_{1} \cup \lbrace a_{2} \rbrace$, $\Omega_{3}= \Omega_{2} \cup \lbrace a_{3} \rbrace \ldots$. By Riemann's removable singularity theorem and maximum modulus principle, we can find holomorphic functions $F_{i}: \Omega_{i} \to \mathbb{D}$ with $F_{i}(\zeta) = f(\zeta)$ for all $\zeta\in \Omega$ and $i\in \mathbb{N}$. Since $\Omega_{i}\subseteq \mathbb{D}$ is a sequence of domains such that $\cup_{i}\Omega_{i} = \mathbb{D}$ and $\Omega_{i} \subseteq \Omega_{i+1}$ for all $i\in \mathbb{N}$, by Montel's theorem, we may assume the sequence $F_{i}$ converges uniformly on each compact subset of $\mathbb{D}$ to a holomorphic function $F:\mathbb{D}\to \overline{\mathbb{D}}$.

Observe that $F(\zeta) = f(\zeta)$ for all $\zeta\in \Omega$ and by maximum modulus principle $F:\mathbb{D}\to \mathbb{D}$. It is easy to see that $F(w)\notin f(\Omega)$ for all $w\in A$. This implies
\begin{align}
S_{\Omega}(z) = c \leq |F(w)|
\end{align}
for all $w\in A$. Now choose a conformal map $\phi : \mathbb{D} \to \mathbb{D}$ such that
$$\phi(\zeta) = \dfrac{\zeta-z}{1-\overline{z}\zeta}.$$
Since Poincaré metric is invariant under conformal self map of $\mathbb{D}$
\begin{align*}
d_{\rho}(z , w)= d_{\rho}(\phi(z) , \phi(w)) = d_{\rho}\left(0 , \dfrac{w-z}{1-\overline{z}w}\right).
\end{align*}
By distance decreasing property of Poincaré metric
\begin{align*}
d_{\rho}(F(z) , F(w))&\leq d_{\rho}(z , w) = d_{\rho}\left(0 , \dfrac{w-z}{1-\overline{z}w}\right)\\
d_{\rho}(0 , F(w))&\leq d_{\rho}\left(0 , \dfrac{w-z}{1-\overline{z}w}\right)\\
\sigma(|F(w)|)&\leq \sigma\left(\left|\dfrac{w-z}{1-\overline{z}w}\right|\right).
\end{align*}
Therefore
\begin{align}
|F(w)|\leq \left|\dfrac{w-z}{1-\overline{z}w}\right|
\end{align}
for all $w\in A$. By equations 2.2 and 2.3
\begin{align}
S_{\Omega}(z) \leq \inf_{w\in A} \left\lbrace \left|\dfrac{w-z}{1-\overline{z}w}\right|\right\rbrace.
\end{align}
Hence the result follows from equations 2.1 and 2.4.
\end{proof}

\begin{remark}
In \cite[Theorem~5.8]{2012}, it is known that for a finitely connected domain $\Omega \subseteq \mathbb{C}$, if any one component of the complement $\Omega$ in $\overline{\mathbb{C}}$ is a single point, then $\Omega$ is not HHR domain. Theorem 1.1 implies that the infinitely connected domain $\Omega = \mathbb{D}\setminus A$ is not HHR domain.  
\end{remark}

\begin{proof}[\underline{Proof of Theorem \rm{\ref{main2}}}]
By {\cite[Proposition~1]{kaushal}} and Theorem 1.1, it is clear that
\begin{align}
h_{\Omega}^{c}(z) \geq \inf_{w\in A} \left\lbrace \left|\dfrac{w-z}{1-\overline{z}w}\right|\right\rbrace.
\end{align}
For $z\in \Omega$, let $f: \mathbb{D}\to \Omega$ be any holomorphic embedding with $\mathbb{B}_{\Omega}^{c}(z, r)\subseteq f(\mathbb{D}) \subseteq \Omega$. This implies that
\begin{align}
\lbrace \zeta\in \Omega : c_{\Omega}(z, \zeta) < r\rbrace \subseteq \Omega.
\end{align}
Now by taking inclusion map from $\Omega$ to $\mathbb{D}$, it is clear that $c_{\mathbb{D}}(\zeta_{1}, \zeta_{2})\leq c_{\Omega}(\zeta_{1}, \zeta_{2})$ for all $\zeta_{1}, \zeta_{2}\in \Omega$. For any holomorphic function $f: \Omega\to \mathbb{D}$, by Riemann's removable singularity theorem and Montel's theorem, we have a holomorphic function $F:\mathbb{D}\to \mathbb{D}$ with $F(\zeta) = f(\zeta)$ for all $\zeta\in \Omega$(similar to Theorem 1.1). This implies $c_{\mathbb{D}}(\zeta_{1}, \zeta_{2})\geq c_{\Omega}(\zeta_{1}, \zeta_{2})$ for all $\zeta_{1}, \zeta_{2}\in \Omega$. Hence $c_{\mathbb{D}}(\zeta_{1}, \zeta_{2})
= c_{\Omega}(\zeta_{1}, \zeta_{2})$ for all $\zeta_{1}, \zeta_{2}\in \Omega$. By equation $2.6$, $\lbrace \zeta\in \Omega : c_{\mathbb{D}}(z, \zeta) < r\rbrace \subseteq \Omega$. Therefore
\begin{align}
\left\lbrace \zeta\in \Omega : \left |\dfrac{\zeta-z}{1-\overline{z}\zeta}\right|< \tanh r \right\rbrace \subseteq \Omega.
\end{align}
Choose a conformal map $\phi :\mathbb{D} \to \mathbb{D}$ such that 
$$\phi(\zeta) = \dfrac{\zeta-z}{1-\overline{z}\zeta}.$$
By equation 2.7, $\lbrace \zeta\in \mathbb{D} : |\zeta| < \tanh r\rbrace \subseteq \phi(\Omega)$. This implies $\tanh r \leq \inf_{w\in A} \left\lbrace \left|\dfrac{w-z}{1-\overline{z}w}\right|\right\rbrace$. Hence the result follows.
\end{proof}
\begin{remark} 
Thus for infinitely connected planar domain $\Omega = \mathbb{D}\setminus A$, $S_{\Omega}(z) = h_{\Omega}^{c}(z)$. Note that in Theorem 1.1, we are taking $A$ as a sequence in $\mathbb{D}$ which converges to the boundary of $\mathbb{D}$. One can see that the formula holds by taking $A$ as a sequence in $\mathbb{D}$ which has finitely many limit points on the boundary of $\mathbb{D}$. It will be interesting to see, what will happen if sequence converges inside $\mathbb{D}$. 
\end{remark}

\begin{remark}
Proof of Theorem 1.1 and Theorem 1.4 can be adopted to show that the squeezing function and Fridman invariant on $\Omega$, where  $\Omega = \mathbb{D}\setminus \{a_{1}\}\cup \{a_{2}\}\cup\ldots\cup\{a_{n}\}$ and $\lbrace a_{i} \rbrace\in \mathbb{D}$ for each $i=1, 2,\ldots, n$, are given by
$$S_{\Omega}(z) = h_{\Omega}^{c}(z) =\min \left\lbrace \left|\dfrac{a_{1}-z}{1-\overline{z}a_{1}}\right|,\ldots, \left|\dfrac{a_{n}-z}{1-\overline{z}a_{n}}\right|\right\rbrace.$$
\end{remark}

\section{Further examples and observations}
\begin{example}\label{ex: 1}
Let $\Omega = \mathbb{D}^{n}\setminus A$ be a domain, where $A = \lbrace a_{k} = (a_{1}^{k}, a_{2}^{k},\ldots, a_{n}^{k}) : k\in \mathbb{N} \rbrace $ is a sequence in $\mathbb{D}^{n}$ converging to $\partial\mathbb{D}^{n}$. Then for $z =(z_{1}, z_{2},\ldots, z_{n})\in \Omega$, we have
$$T_{\Omega}(z) = \inf_{k\in \mathbb{N}} \left\lbrace \max_{1 \leq j\leq n}\left|\dfrac{a_{j}^{k}-z_{j}}{1-\overline{z_{j}}a_{j}^{k}}\right| \right\rbrace.$$
By taking automorphism of $\mathbb{D}^{n}$, it is obvious that 
\begin{align}
T_{\Omega}(z) \geq \inf_{k\in \mathbb{N}} \left\lbrace \max_{1 \leq j\leq n}\left|\dfrac{a_{j}^{k}-z_{j}}{1-\overline{z_{j}}a_{j}^{k}}\right| \right\rbrace.
\end{align}
For $z=(z_{1}, z_{2},\ldots, z_{n})\in \Omega$, let $T_{\Omega}(z) = c$ for some $c\in (0 ,1]$, then using \cite[Theorem~2]{polydisc}, there exists a holomorphic embedding $f=(f_{1}, f_{2},\ldots,f_{n}) :\Omega\to \mathbb{D}^n$ such that $f(z)=0$ with $\mathbb{D}^n(0,c)\subseteq f(\Omega)$. By using Riemann's removable singularity theorem, Montel's theorem and maximum principle on the components of $f$, we get a holomorphic map $F=(F_{1}, F_{2},\ldots,F_{n}):\mathbb{D}^{n}\to \mathbb{D}^n$ such that $F(\zeta)=f(\zeta)$ for all $\zeta \in \Omega$.  

It is easy to see that $F(A)\cap f(\Omega)= \emptyset$. This implies for each $a_{k} \in A$, there is some $F_{i}$ such that $c \leq |F_{i}(a_{k})|$ for $1 \leq i \leq n$. By distance decreasing property of Carathéodory metric
\begin{align*}
c_{\mathbb{D}^{n}}(F(z) , F(a_{k}))&\leq c_{\mathbb{D}^{n}}(z , a_{k})\\
\tanh [c_{\mathbb{D}^{n}}(0 , F(a_{k}))] &\leq \tanh [c_{\mathbb{D}^{n}}(z , a_{k})]\\
\max_{1\leq i \leq n}|F_{i}(a_{k})|&\leq \tanh [c_{\mathbb{D}^{n}}(z , a_{k})].
\end{align*}
Therefore
\begin{align}
T_{\Omega}(z) \leq \inf_{k\in \mathbb{N}} \lbrace \tanh [c_{\mathbb{D}^{n}}(z , a_{k})]\rbrace.
\end{align}
Hence the result follows from equations 3.1 and 3.2.
\end{example}

\begin{remark}
Similarly one can see that for domain $\Omega = \mathbb{D}^{n}\setminus \{a_{1}\}\cup \{a_{2}\}\cup\ldots\cup\{a_{k}\}$, where $\lbrace a_{i} \rbrace\in \mathbb{D}^{n}$ for each $i=1, 2,\ldots, k$(here $a_i=(a_{i1}, a_{i2},\ldots, a_{in})$) and for $z =(z_{1}, z_{2},\ldots, z_{n})\in\Omega$, we have
$$T_{\Omega}(z) = \min \left\lbrace \max_{1 \leq j\leq n}\left|\dfrac{a_{1j}-z_{j}}{1-\overline{z_{j}}a_{1j}}\right|,\ldots, \max_{1 \leq j\leq n} \left|\dfrac{a_{kj}-z_{j}}{1-\overline{z_{j}}a_{kj}}\right|\right\rbrace.$$
Although we can get this by {\cite[Theorem~2.1]{rong3}}.
\end{remark}

\begin{lemma}{\cite[Theorem~1.2.6]{riemann-removable}}(Hartogs’s extension theorem)
Let $\Omega$ be a domain in $\mathbb{C}^{n}(n \geq 2)$. Suppose that $A$ is a compact subset of $\Omega$ with $\Omega\setminus A$ is connected. If $f$ is holomorphic on $\Omega\setminus A$, then there is a holomorphic function $F$ on $\Omega$ such that $F_{\Omega\setminus A}=f$.
\end{lemma}

\begin{example}\label{ex: 2}
Let $\lbrace a_{k} : k\in \mathbb{N}\rbrace$ be a sequence in $\mathbb{D}^{n}(n\geq 2)$ converging to $\partial\mathbb{D}^{n}$. Suppose that $\mathbb{D}_{k}^{n} \subset \mathbb{D}^{n}$ is a polydisk centred at $a_{k}$ with radius $r_{k}$ for each $k\in \mathbb{N}$ such that $\mathbb{D}_{i}^{n}\cap \mathbb{D}_{j}^{n} = \emptyset$ for any $i, j \in \mathbb{N}$. Let $\Omega = \mathbb{D}^{n}\setminus A$ be a domain, where $A = \cup_{k\in \mathbb{N}}\overline{\mathbb{D}_{k}^{n}}$. Then for $z =(z_{1}, z_{2},\ldots, z_{n})\in \Omega$, we have
$$T_{\Omega}(z) = \inf_{k\in \mathbb{N}} \left\lbrace \max_{1 \leq j\leq n} \left|\dfrac{r_{k}-|z_{j}|}{1-|z_{j}|r_{k}}\right| \right\rbrace.$$
For $z =(z_{1}, z_{2},\ldots, z_{n})\in \Omega$, let $f\in Aut(\mathbb{D}^{n})$ with $f(z) = 0$ and $ \lambda = \inf_{k \in \mathbb{N}} \left\lbrace \min_{w\in \partial\mathbb{D}_{k}^{n}}\tanh [c_{\mathbb{D}^{n}}(z , w)]\right\rbrace$. Clearly $g = f|_{\Omega}$ is a holomorphic embedding from $\Omega$ to $\mathbb{D}^{n}$ with $g(z) = 0$. We claim that $\mathbb{D}^{n}(0, \lambda)\subseteq g(\Omega)$. To prove this, let $\alpha \in \mathbb{D}^{n}(0, \lambda)$. This implies  $\tanh c_{\mathbb{D}^{n}}(0 , \alpha) < \lambda$ and since $f$ is an automorphism, therefore $\tanh c_{\mathbb{D}^{n}}(0 , \alpha) = \tanh c_{f(\mathbb{D}^{n})}(f(z) , f(\alpha')) = \tanh c_{\mathbb{D}^{n}}(z , \alpha')$. 

Thus $\tanh c_{\mathbb{D}^{n}}(z , \alpha') < \inf_{k \in \mathbb{N}} \left\lbrace \min_{w\in \partial\mathbb{D}_{k}^{n}}\tanh [c_{\mathbb{D}^{n}}(z , w)]\right\rbrace$. This shows $\alpha'\notin A$ and hence $\alpha = f(\alpha')$ for some $\alpha'\in \Omega$. This proves our claim. Therefore
\begin{align}
T_{\Omega}(z) \geq \inf_{k \in \mathbb{N}} \left\lbrace \min_{w\in \partial\mathbb{D}_{k}^{n}}\tanh [c_{\mathbb{D}^{n}}(z , w)]\right\rbrace.
\end{align}
By replacing Riemann's removable singularity theorem with Hartogs’s extension theorem in Example 3.1, it is easy to see that
\begin{align}
T_{\Omega}(z) \leq \inf_{k \in \mathbb{N}} \left\lbrace \min_{w\in \partial\mathbb{D}_{k}^{n}}\tanh [c_{\mathbb{D}^{n}}(z , w)]\right\rbrace.
\end{align}
By equations 3.3 and 3.4, we have
\begin{align*}
T_{\Omega}(z) = \inf_{k\in \mathbb{N}} \left\lbrace \max_{1 \leq j\leq n} \left|\dfrac{r_{k}-|z_{j}|}{1-|z_{j}|r_{k}}\right| \right\rbrace.
\end{align*}
\end{example}

\begin{example}\label{ex: 3}
Let $\lbrace a_{k} : k\in \mathbb{N}\rbrace$ be a sequence in $\mathbb{D}^{n}(n\geq 2)$ converging to $\partial\mathbb{D}^{n}$. Suppose that $\mathbb{B}_{k}^{n} \subset \mathbb{D}^{n}$ is a ball centred at $a_{k}$ with radius $r_{k}$ for each $k\in \mathbb{N}$ such that $\mathbb{B}_{i}^{n} \cap \mathbb{B}_{j}^{n} = \emptyset$ for any $i, j \in \mathbb{N}$. Let $\Omega = \mathbb{D}^{n}\setminus A$ be a domain, where $A = \cup_{k\in \mathbb{N}}\overline{\mathbb{B}_{k}^{n}}$. Then for $z =(z_{1}, z_{2},\ldots, z_{n})\in \Omega$, we have
$$T_{\Omega}(z) = \inf_{k\in \mathbb{N}} \left\lbrace \min_{w\in \partial\mathbb{B}_{k}^{n}}\max_{1 \leq j\leq n} \left|\dfrac{|w_{j}|-|z_{j}|}{1-|z_{j}||w_{j}|}\right| \right\rbrace.$$
Very similar to Example 3.4, we omit the details.
\end{example}

In \cite{polydisc}, there is a statement which tells that for bounded homogeneous domain $\Omega\subseteq\mathbb{C}^{n}$, 
either $S_{\Omega}(z)= \dfrac{1}{\sqrt{n}}T_{\Omega}(z)$ or $T_{\Omega}(z)= \dfrac{1}{\sqrt{n}}S_{\Omega}(z)$. 
We see that this statement is not true in general. To see this consider the following Remark. 
\begin{remark} 
Let $\Omega = \mathbb{B}^{n}\times\mathbb{B}^{n}\times\ldots\times\mathbb{B}^{n}\subset \mathbb{C}^{n^{2}}$ be a domain.
 Then for $z=(z_{1}, z_{2},\ldots, z_{n})\in \Omega$ and $n > 1$(here $z_i=(z_{i1}, z_{i2},\ldots, z_{in})$), we have
$$S_{\Omega}(z)\neq \dfrac{1}{n}T_{\Omega}(z)\  \ \ \mbox{and}\  \ \ T_{\Omega}(z)\neq\dfrac{1}{n} S_{\Omega}(z).$$
Since $\mathbb{B}^{n}$ is a classical symmetric domain, by \cite[Theorem~7.5]{2012}
$$S_{\Omega}(z) = (S_{\mathbb{B}^{n}}(z_{1})^{-2}+S_{\mathbb{B}^{n}}(z_{2})^{-2}+\ldots +S_{\mathbb{B}^{n}}(z_{n})^{-2})^{-1/2} = \dfrac{1}{\sqrt{n}}.$$
Also, by \cite[Proposition~1]{polydisc}, $T_{\Omega}(z)\geq\displaystyle\min_{1\leq i \leq n}T_{\mathbb{B}^{n}}(z_{i})$ and since $T_{\mathbb{B}^{n}}(z_{i})\equiv \dfrac{1}{\sqrt{n}}$ for each $i$, this implies $T_{\Omega}(z)\geq \dfrac{1}{\sqrt{n}}$.\\ 
If $S_{\Omega}(z)= \dfrac{1}{n}T_{\Omega}(z)$, then $T_{\Omega}(z)=\sqrt{n}$ which is not possible.\\
If $T_{\Omega}(z)= \dfrac{1}{n}S_{\Omega}(z)$, then $T_{\Omega}(z)=\dfrac{1}{n^{3/2}} < \dfrac{1}{\sqrt{n}}$ which is not possible. \\
Note that $\Omega = \mathbb{B}^{n}\times\mathbb{B}^{n}\times\ldots\times\mathbb{B}^{n}\subset \mathbb{C}^{n^{2}}$ is a bounded homogeneous domain.
\end{remark}

In {\cite[Theorem~2.1]{rong3}}, Rong and Yang have given an explicit formula of squeezing function for domain $\Omega\setminus K$, where $\Omega$ is a bounded, balanced, convex and homogeneous domain in $\mathbb{C}^{n}$ and $K$ is a compact subset of $\Omega$. Although they have not mentioned, for the sake of completion we have given an example to show that the formula will not hold in case of planar domains. 
\begin{remark}
Let $\Omega = \mathbb{D}$ and $K = \overline{\mathbb{D}}_{1/4} = \left\lbrace z\in \mathbb{C} : |z|\leq \dfrac{1}{4} \right\rbrace$, then $\Omega\setminus K = A_{1/4} = \left\lbrace z\in \mathbb{C} : \dfrac{1}{4}<|z|<1 \right\rbrace$. For $z=\dfrac{1}{2}$, by {\cite[Theorem~2.1]{rong3}} 
\begin{align*}
S_{A_{1/4}}\left(\dfrac{1}{2}\right) = \min_{w\in \overline{\mathbb{D}}_{1/4}} \tanh\left[c_{\mathbb{D}}\left(\dfrac{1}{2}, w\right)\right]
\end{align*}
This implies
\begin{align*}
S_{A_{1/4}}\left(\dfrac{1}{2}\right) = \min_{w\in \overline{\mathbb{D}}_{1/4}} \left|\dfrac{w-\dfrac{1}{2}}{1- \dfrac{w}{2}}\right|.
\end{align*}
By choosing $w=\dfrac{1}{4}$, we have $S_{A_{1/4}}\left(\dfrac{1}{2}\right) \leq \dfrac{2}{7}$ but by the expression of squeezing function on annulus, $S_{A_{1/4}}\left(\dfrac{1}{2}\right) = \dfrac{1}{2}$.
\end{remark}

\begin{remark}
It can be seen that Theorem 1.1 and Theorem 1.4 are special case of Theorem 2.8 of \cite{rong3}. We would like to clarify here that we wanted to give an explicit and easy proof without the technicality of analytic sets in planar domains. We are thankful to Shichao Yang for reminding these things.
\end{remark}

\section*{Ackowledgement}
We are thankful to our thesis advisor Sanjay Kumar and co-researcher Naveen Gupta for his valuable suggestions.

\end{document}